\pdfoutput=1
\RequirePackage{ifpdf}
\ifpdf 
\documentclass[pdftex]{sigma}
\else
\documentclass{sigma}
\fi

\begin{document}

\renewcommand{\thefootnote}{}

\newcommand{\arXivNumber}{2501.18061}

\renewcommand{\PaperNumber}{015}

\FirstPageHeading

\ShortArticleName{Experimenting with the Garsia--Milne Involution Principle}

\ArticleName{Experimenting with the Garsia--Milne Involution\\ Principle\footnote{This paper is a~contribution to the Special Issue on Basic Hypergeometric Series Associated with Root Systems and Applications in honor of Stephen C.~Milne's 75th birthday. The~full collection is available at \href{https://www.emis.de/journals/SIGMA/Milne.html}{https://www.emis.de/journals/SIGMA/Milne.html}}}

\Author{Shalosh B.~EKHAD and Doron ZEILBERGER}

\AuthorNameForHeading{S.B.~Ekhad and D.~Zeilberger}

\Address{Department of Mathematics, Rutgers University (New Brunswick), \\
Hill Center-Busch Campus, 110 Frelinghuysen Rd., Piscataway, NJ 08854-8019, USA}
\Email{\href{mailto:DoronZeil@gmail.com}{DoronZeil@gmail.com}}
\URLaddress{\url{https://sites.math.rutgers.edu/~zeilberg/}}

\ArticleDates{Received January 31, 2025, in final form February 26, 2025; Published online March 04, 2025}

\Abstract{In 1981, Adriano Garsia and Steve Milne found the first bijective proof of the celebrated Rogers--Ramanujan identities. To achieve this feat, they invented a versatile tool that they called the Involution Principle. In this note we revisit this useful principle from a~very general perspective, independent of its application to specific combinatorial identities, and will explore its complexity.}


\Classification{05A19}

\begin{flushright}
\begin{minipage}{90mm}
\it In memory of Adriano M. Garsia (1928--2024)\\ and in honor of Stephen C. Milne on his 75th birthday
\end{minipage}
\end{flushright}

\renewcommand{\thefootnote}{\arabic{footnote}}
\setcounter{footnote}{0}

\section{The ``million dollar'' problem in partition theory back in 1980}

One of the most famous identities in {\it both} {\it enumerative combinatorics} and {\it number theory} are the {\it Rogers--Ramanujan identities} (see Drew Sills' wonderful book \cite{S} for their fascinating history and
for many references).
Already by 1940, G.H. Hardy knew of {\it seven} different proofs. In his own words~\cite{H}:

{\rm ``I can remember very well his surprise, and the admiration which he expressed for Rogers' work. A correspondence followed in the course of which Rogers was led to a considerable simplification of his original proof.
About the same time I.~Schur, who was then cut off from England by the war, rediscovered the identities again. Schur published two proofs, one of which is `combinatorial' and quite unlike any other proof known.
There are now seven published proofs, the four referred to already, the two much simpler proofs found later by Rogers and Ramanujan and published in the papers,
and a much later proof by Watson based on quite different ideas. None of these proofs can be called both ``simple'' and ``straightforward'', since the simplest are essentially verifications; and no doubt it would be unreasonable to expect a really easy proof.''
}

Forty years later, in the early 1980s, there were a few additional proofs, including a lovely combinatorial proof by Dave Bressoud \cite{B}. But none of them was {\it bijective}.
Recall that the first Rogers--Ramanujan identity is equivalent to the fact that the cardinality of the set of integer partitions of $n$ into parts
that differ from each other by at least $2$, let us call it $A(n)$, has the same number of elements as the set of integer partitions of $n$ whose parts are
congruent to one or four modulo five. The great open problem at that time was to find an {\it explicit bijection} from the set $A(n)$ onto the set $B(n)$. This
means to come up with an {\it algorithm}, let us call it~$\phi$, that inputs a member $a \in A(n)$ and outputs a member $b \in B(n)$, and does {\it not} use the fact
that~${|A(n)|=|B(n)|}$, thereby giving a {\it bijective proof} of this identity.

The (already then) {\it grand old man} of partition theory, George Andrews, offered a monetary prize for this challenge, that Adriano Garsia and Stephen Milne
proudly won \cite{GM1, GM2}. Fairly recently, George Andrews, in one of the fascinating interviews conducted by the visionary enumerative enumerator Toufik Mansour,
commented \cite{AM}:

{\rm ``Finally, and most importantly, the $q$-world
is the natural home of the generating functions
for integer partitions. For example, the celebrated Rogers--Ramanujan identities were first
proved in~$1894$ by Rogers. It took until $1981$
when Garsia and Milne published a purely bijective proof in a $50$-page paper.''
}

Very recently (January 24, 2025) George Andrews confirmed it in an email message, that he kindly allowed us to quote:

{\rm
``As I recall, Steve called me and said they had a proof. I was skeptical until I heard the word `non-canonical'. That meant to me that their algorithm would somehow blindly match things up but, at the end, you would have no feeling for
what happened, nor would you have any idea of a refinement of Rogers--Ramanujan. I wrote the check for $50$ dollars on the spot and sent it to Milne that evening. I think this was in the early $1980$.''}

Like Christopher Columbus \cite{Wik1} and Alexander the Great \cite{Wik2} before them, Garsia and Milne had a brilliant idea: {\it think outside the box!}
First come up with two much larger (but still finite) sets $X(n)$ and $Y(n)$ such that there are {\it natural} bijections,
$
\phi\colon  X(n) \rightarrow Y(n) ,
$
and~${
\psi\colon   X(n)\backslash A(n) \rightarrow Y(n)\backslash B(n)}$,
then construct the {\it artificial} bijection
$
\alpha\colon   A(n) \rightarrow B(n) ,
$
as follows. Given $a \in A(n)$, let $k$ be the smallest non-negative integer such that
\smash{$
\bigl(\phi \psi^{-1}\bigr)^k\phi(a) \in B(n)$},
 and define
\smash{$
\alpha(a):=\bigl(\phi \psi^{-1}\bigr)^k\phi(a)$}.

{\bf Comment.} The original set-up for the involution principle (hence the name), succinctly summarized by Jeff Remmel (see \cite[Section 1]{R}), involved a pair of involutions $\alpha$ and $\beta$ and sets~${A=A^{+} \cup A^{-}}$,
$B=B^{+} \cup B^{-}$ as well as a sign-preserving bijection from $A$ to $B$ and {\it sign-reversing} involutions $\alpha$ and $\beta$ defined on $A$ and $B$ respectively. It is easy to see that
this scenario can be reduced to the above simplified version.

Artificial or not it was {\it brilliant} and became {\it iconic}. It was deemed important enough to be included in the chapter ``Enumerative and algebraic combinatorics'' \cite{Z} in
the monumental {\it Princeton companion to mathematics}, edited by Sir Timothy Gowers. Let us quote:

{\rm ``When we perform algebraic $($and logical, and
even analytical$)$ manipulations, we are really rearranging and
{\it combining} symbols, hence doing combinatorics in disguise.
In fact, {\it everything is
combinatorics}. All we need to do is to take the combinatorics out of the
closet, and make it explicit.
The plus sign
turns into $($disjoint$)$ union, the multiplication sign becomes
Cartesian product, and induction turns into recursion. But
what about the combinatorial counterpart of the minus sign?

In $1982$, Adriano Garsia and Stephen Milne filled
this gap by producing an
ingenious `Involution principle' that enables one to
translate the implication
\[
a=b \qquad \text{and} \qquad c=d \qquad \Rightarrow \qquad
a-c=b-d ,
\]
into a bijective argument,
in the sense that if
$C \subset A$ and $D \subset B$, and
there are natural bijections~${f\colon A \rightarrow B}$ and $g\colon C \rightarrow D$ establishing
that $a:=\vert A \vert$ equals $b:=\vert B \vert$, and
$c:=\vert C \vert$ equals~${d:=\vert D \vert}$, then it is
possible to construct an explicit bijection between
$A \backslash C$ and~$B \backslash D $. Let us define it in terms of people.
Suppose that in a certain village all the adults are
married, and hence there is a natural bijection between
the set of married men to the set of married women, ${m \rightarrow \operatorname{WifeOf}(m)}$, with its inverse
$w \rightarrow \operatorname{HusbandOf}(w)$.
In addition, some of the people have extra-marital affairs,
$($but only one per person, and all within the village$)$.
There is a natural bijection between the set of cheating men
to the set of cheating women, called~${m \rightarrow \operatorname{MistressOf}(m)}$,
with its inverse $w \rightarrow \operatorname{LoverOf} (w)$.
It follows that there are as many faithful men as there are
faithful women. How to match them up? $($say if
a faithful man wants a faithful woman to go to Church with him$)$.
A faithful man first asks his wife to come with him.
If she is faithful, she agrees.
If she is not, she has a lover, and that lover has a wife.
So she tells her husband: sorry, hubby, but I am going to the
pub with my lover, but my lover's wife may be free,
so the man asks the wife of the lover of his wife to go with him,
and if she is faithful, she agrees. If she is not,
he keeps asking the wife of the lover of the woman who
just rejected his proposal, and since the village is finite,
he will eventually get to a faithful woman.

The reaction of the combinatorial enumeration community to
the involution principle was mixed. On the one hand it had
the {\it universal appeal} of a general {\it principle}, that
should be useful in many attempts to find bijective
proofs of combinatorial identities. On the other hand,
its universality is also a major drawback, since involution principle
proofs usually do not give any insight into the {\it specific}
structures involved, and one feels a bit cheated.
Such a proof answers
the {\it letter} of the question, but misses its {\it spirit}.
In these cases
one still hopes for a~{\it really} natural, `involution principle-free proof'. This is the case with the celebrated Rogers--Ramanujan
identity that states that the number of partitions of an integer
into parts that leave remainder $1$ or~$4$ when divided by $5$ equals
the number of partitions of that integer with the property
that the difference between
parts is at least $2$. For example, if $n=7$ the
cardinalities of $\{61, 4111, 1111111\}$ and~${\{7,61,52\}}$
are the same.
Garsia and Milne invented their notorious
{\it principle} in order to give a Rogers--Ramanujan bijection,
thereby winning a $50$ dollar prize from George Andrews.
However, finding a {\it really nice} bijective proof is still an
open problem.''}

\section[Implementing and experimenting with random Garsia--Milne scenarios]{Implementing and experimenting\\ with random Garsia--Milne scenarios}

In order to utilize their ingenious new principle, Garsia and Milne \cite{GM1,GM2} had to come up with a~{\it specific} set-up for which it could be applied.
Jeff Remmel \cite{R} also applied it brilliantly to many other scenarios. This was further extended to {\it general} partition identities by Kathy O'Hara~\cite{O},
Herb Wilf~\cite{Wil}, David Feldman and Jim Propp~\cite{FP}, and others. We should also mention the interesting approach of Ilse Fischer and Matja\v{z} Konvalinka~\cite{FK},
who turned the involution principle into a calculus of `signed sets' and `sijections' (as they call them).

Here we will {\it forget} about the origin of the involution principle
in partition theory and study it completely generally, using the very general framework of cheating and faithful men and women, where
the assignment of the mapping ``mistress of" is completely {\it arbitrary} and random. Suppose you have $c$ cheating men and $f$ faithful men,
and also $c$ cheating women and $f$ faithful women.
Without loss of generality the cheating men are Mr.~$i$, for $1\leq i \leq c$ and the faithful men are Mr.~$i$ for $c+1 \leq i \leq c+f$.
Both marital and extra-marital relationships are monogamous. So any such situation can be described as a list of length $c$
$[M_1, \dots, M_c]$,
where Mrs.~$M_i$ is the mistress of Mr.~$i$. The entries are all distinct and belong to $\{1, \dots, c+f\}$.

Note that there are $(c+f)(c+f-1) \cdots (f+1)$ ways of doing it.

We were wondering about the {\it complexity} (i.e., the number of requests) of the resulting mapping from faithful men to faithful women, both for
one specific individual (without loss of generality, Mr.~$c+1$) and also about the total complexity adding up all the requests by all faithful men.

\section[The statistical distribution of the number of requests of ONE specific faithful man]{The statistical distribution of the number of requests\\ of ONE specific faithful man}

For any specific faithful man (w.l.o.g.~Mr.~$c+1$) what is the
probability that he would have to make $i$ requests? If his wife is
faithful (has no lover) it is one request. The other extreme is~${i=c+1}$, where he had
to ask all the cheating women, $c$ of them, until he got a match. It
turns out (see below) that this probability is
\[
\frac{{{f+c-i} \choose {f-1}}}{{{c+f} \choose {c}}} .
\]
This is a {\it hypergeometric} distribution whose {\it mean} is
$\frac{c +f +1}{f+1}$.

Note that if the number of {\it cheating} men far exceeds the number of {\it faithful} men (that is definitely the case in all the applications to partition identities) this
average complexity is rather large and tends to the ratio $c/f$.

The variance is
\[
\frac{(c +f +1) c f}{(f +1)^{2} (2+f )} .
\]

The third central moment is
\[
\frac{(c +f +1) c f \big(2 c f +f^{2}-2 c -1\big)}{(f +1)^{3} (2+f ) (3+f )} .
\]

The fourth central moment is\footnote{For higher central moments see the output file
\url{https://sites.math.rutgers.edu/~zeilberg/tokhniot/oGMIP2s.txt}.}
\[
\frac{(c +f +1) c f \big(9 c^{2} f^{2}+9 c  f^{3}+f^{4}-3 c^{2} f +6 c  f^{2}-3 f^{3}+6 c^{2}+3 c f -9 f^{2}+6 c -5 f \big)}{(f +1)^{4} (2+f ) (3+f ) (4+f )} .
\]

Letting $c=f$ it is easily seen that the {\it limiting distribution} as $f$ goes to infinity is the Geometric distribution $\operatorname{Ge}(1/2)$. More generally
for any fixed $k$ setting $c=kf$ and letting $f$ go to infinity, the limiting distribution is $\operatorname{Ge}(1/(k+1))$.

\section[The statistical distribution of the number of total requests by all faithful men]{The statistical distribution of the number of total requests\\ by all faithful men}

What about the {\it total} number of requests by all faithful men until they find a woman to go to church with?

This number can be anything between $f$ (where all the faithful men are married to faithful women) and $c+f$. The probability that it is $i$ equals
\[
\frac{{{i-1} \choose {f-1}}}{{{c+f} \choose {c}}} .
\]
The {\it mean} is, by {\it linearity of expectation}, $f$ times the mean for one individual, given above, hence
$
\frac{f (c +f +1)}{f +1} .
$

The variance is
\[
\frac{(c +f +1) c f}{(f +1)^{2} (2+f )} .
\]

The third central moment is
\[
-\frac{(c +f +1) c f \big(2 c f +f^{2}-2 c -1\big)}{(f +1)^{3} (2+f ) (3+f )} .
\]

The fourth central moment is\footnote{For higher central moments see the output file
\url{https://sites.math.rutgers.edu/~zeilberg/tokhniot/oGMIP1.txt}.}
\[
\frac{(c +f +1) c f \big(9 c^{2} f^{2}+9 c  f^{3}+f^{4}-3 c^{2} f +6 c  f^{2}-3 f^{3}+6 c^{2}+3 c f -9 f^{2}+6 c -5 f \big)}{(f +1)^{4} (2+f ) (3+f ) (4+f )} .
\]

\section{A brief explanation}

All the above results were obtained via the Maple package {\tt GMIP.txt} written by the second author and executed by the first one.\footnote{It is available from \url{https://sites.math.rutgers.edu/~zeilberg/tokhniot/GMIP.txt}.}

Denoting the lengths of the paths of the $f$ faithful men until they reach their matched faithful woman by
$
[a_1 , \dots , a_f] ,
$
it is readily seen that
$
f \leq a_1+ \dots + a_f \leq c+f .
$

Each such assignment of path-lengths can occur in $f!$ ways. Hence assigning the weight $t_{c+1}^{a_1} \cdots t_{c+f}^{a_f}$
to each matching, the {\it weight-enumerator} is the coefficient of $z^c$ in the Taylor expansion of
\[
\frac{f!}{1-z} \cdot \prod_{i=c+1}^{c+f} \frac{t_{i}}{1-z t_{i}} .
\]
If we want to focus on one specific man, say Mr. $c+1$, we set $t_{c+1}=x$ and $t_{c+2}= \dots =t_{c+f}=1$, and {\it normalize}, in order to convert it into a {\it probability generating function}, we get that
it equals the coefficient of $z^c$ in the Taylor expansion of
\[
\frac{c!f!}{(c+f)!} \cdot \frac{x}{(1-xz)(1-z)^f} .
\]

If we want to focus on the total complexity, you set all the $t_i=x$ getting that the {\it probability generating function} is the coefficient of $z^c$ in the Taylor expansion of
\[
\frac{c!f!}{(c+f)!} \cdot \frac{x^f}{(1-z) (1-xz)^f} .
\]

\section{Encore}

Our package can also spell out any random scenario. Procedure {\tt RandL(n,k)} gives a random assignment when there are $n$ men and $n$ women and $k$ cheating men and $k$ cheating women.
Procedure {\tt GMstory(n,L)} spells out all the requests by the faithful men and gives the resulting matching.
Here is a small example:
Typing {\tt RandL(4,3);} may yield (out of $4 \cdot 3 \cdot 2=24$ possibilities).
$
L=[2,4,3] .
$
Then typing {\tt GMstory(4,L);} tells you the story:

Once upon a time there was a village with $4$ married couples. Let us call them Mr.~$i$ and Mrs.~$i$, where Mr.~$i$ and Mrs.~$i$ are married to each other, where $i$ goes from $1$ to $4$.

 It so happened that Mr.~$1$, Mr.~$2$, and Mr.~$3$ are cheaters. They each have one mistress:
\begin{itemize}\itemsep=0pt
\item The Mistress of Mr.~$1$ is Mrs.~$2$ (and hence the Lover of Mrs.~$2$ is Mr.~$1$)

\item The Mistress of Mr.~$2$ is Mrs.~$4$ (and hence the Lover of Mrs.~$4$ is Mr.~$2$)

\item The Mistress of Mr.~$3$ is Mrs.~$3$ (and hence the Lover of Mrs.~$3$ is Mr.~$3$) (this can happen!)
\end{itemize}
 The $3$ cheating men and their mistresses refuse to go to Church on Sunday, they would rather go to the pub.

 There must be a way to match the only faithful man, Mr.~$4$ with the only faithful woman Mrs.~$1$.

Mr.~$4$ first asks his wife: ``My dear wife, will you go with me to Church?''

She replies: ``Sorry, hubby, but I am going to the pub with my lover, Mr.~$2$, perhaps his wife Mrs.~$2$ is willing?

So Mr.~$4$ asks Mrs.~$2$ and she replies: ``Sorry Mr.~$4$, I can't make it, I am going to the pub with my lover Mr.~$1$,
why won't you ask his wife, Mrs.~$1$?

So he asks Mrs.~$1$, and she happily accepts, since she has no lover and would rather go to Church.
So Mr.~$4$ goes to Church with Mrs.~$1$.

{\bf Simulation.} We also have simulation procedures, that empirically confirm the theoretical results. Enjoy!

\section{Conclusion}

When Adriano Garsia and Stephen Milne astounded the enumerative combinatorics community with the first proof of the Rogers--Ramanujan identities,
in some sense they {\it ``cheated''}, but in a {\it good way}! By doing it they invented an essential tool for constructing bijections
between combinatorial objects that often do not seem to have a natural {\it canonical} mapping between them, so
{\it faute de mieux}, one has to resort to it. Since it is so general, one can forget about the specific combinatorial
framework and play around with random men and women and study the statistical distribution of their interactions.

\subsection*{Acknowledgments} Many thanks to Svante Janson for helpful probability guidance. Also many thanks to anonymous referees
for corrections and insightful comments.

\pdfbookmark[1]{References}{ref}
\LastPageEnding

\end{document}